\documentstyle[amstex,amssymb]{amsart}

\newcommand{\cc}{{\mathbb C}}
\newcommand{\lb}{\lambda}
\newcommand{\td}{{\mathbb T}}

\newcommand{\z}{{\mathbb Z}}
\newcommand{\pn}{{\mathbb P}^n}

\newcommand{{\rr}}{\mathcal{R}}
\newcommand{\bx}{\hfill \rule{.25cm}{.25cm} \medbreak}

\newtheorem{thm}{Theorem}[section]

\newtheorem{corr}[thm]{Corollary}
\newtheorem{prop}[thm]{Proposition}

\begin{document}
\Large

\title{Projective spectrum in Banach algebras}

\author[R. Yang]{Rongwei Yang* }
\address{Department of Mathematics and Statistics, SUNY at Albany, 
Albany, NY 12222, U.S.A.}
\email{ryang@@math.albany.edu}
\thanks{This work is supported in part by a grant from the National Science Foundation (DMS 0500333).} 


\begin{abstract}
\large
For a tuple $A=(A_0,\ A_1,\ ...,\ A_n)$ of elements in a unital Banach algebra ${\mathcal B}$, its {\em projective spectrum} $p(A)$ 
is defined to be the collection of $z=[z_0,\ z_1,\ ...,\ z_n]\in \pn$ such that
$A(z)=z_0A_0+z_1A_1+\cdots +z_nA_n$ is not invertible in ${\mathcal B}$. The pre-image of $p(A)$ in ${\cc}^{n+1}$ is denoted by $P(A)$.
When ${\mathcal B}$ is the $k\times k$ matrix algebra $M_k(\cc)$, the projective spectrum is a projective hypersurface.  In infinite dimensional cases, 
projective spectrums
can be very complicated, but also have some properties similar to that of hypersurfaces. When $A$ is commutative, $P(A)$ 
is a union of hyperplanes. When ${\mathcal B}$ is reflexive or is a $C^*$-algebra,
the {\em projective resolvent set} $P^c(A):=\cc^{n+1}\setminus P(A)$ is shown to be a disjoint union of domains of 
holomorphy. Later part of this paper studies
Maurer-Cartan type ${\mathcal B}$-valued 1-form $A^{-1}(z)dA(z)$ on $P^c(A)$. As a consequence, we show that if ${\mathcal B}$ is a $C^*$-algebra with 
a trace $\phi$, then $\phi(A^{-1}(z)dA(z))$ is a nontrivial element in the de Rham cohomology space $H^1_d(P^c(A),\ \cc)$.
 \end{abstract}

\maketitle

\section{Introduction}

The classical spectrum of an element $A$ in a unital Banach algebra ${\mathcal B}$ is defined through the invertibility of 
$A-\lb I$. If $A=(A_0,\ A_1,\ ...,\ A_n)$ is a commutative tuple of elements in ${\mathcal B}$, then classical notions of joint spectrum 
are defined through the invertibility of $(A_0-\lb_0I,\ A_1-\lb_1I,\ ...,\ A_n-\lb_nI)$ in various senses (H\"{o}rmander [H\"{o}] Ch3, and Taylor [Ta]). 
In all these cases, the identity 
$I$ serves as a base against which the invertibilities of other elements are measured. The idea of projective spectrum, which we will define and 
study, is to set $I$ free, and consider the invertibility of  $z_0A_0+z_1A_1$, or more generally, 
$A(z):=z_0A_0+z_1A_1+\cdots +z_nA_n.$
This is a measurement of how the elements behave against each other. Unlike classical notions of joint spectrums,
projective spectrum is defined for all tuples, not just commutative ones. This paper is organized as follows.\\

Section 1. {\em Preparation}. Here we define the projective spectrum and prove its non-triviality.

Section 2. {\em Projective spectrum and hypersurfaces}. When the Banach algebra is the matrix algebra $M_k(\cc)$, projective spectrums are degree $k$ 
projective hypersurfaces. A comparision
between general projective spectrums and hypersurfaces is made in Section 2. We will see that when $A$ is a commutative tuple, its projective 
spectrum is a union of hyperplanes. The main results are regarding the complement of projective spectrum (which we call projective 
resolvent set). We show that when the Banach algebra is of certain type, for instance $C^*$, the complement is made of domains 
of holomorphy.

Section 3. {\em ${\mathcal B}$-valued 1-form $A^{-1}(z)dA(z)$ and the de Rham cohomology space $H^1_d(P^c(A),\ \cc)$}.
This section makes a study on the topology of projective resolvent set. Since the tuple $A$ in general is of infinite dimensional nature, its 
projective resolvent can be very complicated. Nonetheless, with the aid of Maurer-Cartan type form $A^{-1}(z)dA(z)$ and 
central linear functionals on ${\mathcal B}$, we manage to peek into the de Rham cohomology space of projective resolvent sets.

Section 4. {\em The case when when $A$ is commutative}. A few observations and remarks are made here 
concerning an Arnold and Briskorn's theorem in Hyperplane Arrangements.\\
 
{\bf Acknowlegement:} This paper was benefited from conversations with many of the authors colleagues at SUNY at Albany,
to whom the author is deeply indebted. In particular, the author would like to thank Michael Range, Mark Steinberger,
Michael Stessin, and Alex Tchernev for references and valuable discussions.

\section{Preparation}

We let $z=(z_0,\ z_1,\ ...,\ z_n)$ denote a general point in ${\cc}^{n+1}$. The group ${\cc}^{\times}$ of nonzero complex 
numbers acts on ${\cc}^{n+1}$ by scalar multiplications. The $n$ dimensional projective space $\pn$ is the quotient 
$({\cc}^{n+1}\setminus \{0\})/{\cc}^{\times}$. With topology induced from this quotient, $\pn$ is a compact complex manifold. 
The fibres of this quotient map are the integral curves of the {\em Euler vector field} 
$\theta =\sum_{j=0}^{n}z_j\frac{\partial}{\partial z_j}$. 
For a subset $S\subset \cc^{n+1}$ invariant under ${\cc}^{\times}$, $(S\setminus \{0\})/{\cc}^{\times}$ will be denoted by $S^T$.  
$[z_0,\ z_1,\ ...,\ z_n]$ denotes the homogenous coordinate of a general point in $\pn$. On the open subset $U_0=\{z_0\neq 0\}\subset \pn$, 
$[z_0,\ z_1,\ ...,\ z_n]=[1,\ z_1/z_0,\ z_2/z_0,\ ...,\ z_n/z_0]$. The tuple $(z_1/z_0,\ z_2/z_0,\ ...,\ z_n/z_0)$ is the affine 
coordinate for $U_0$, and is denoted simply by $(\xi_1,\ \xi_2,\ ...,\ \xi_n)$.  

Throughout this paper, ${\mathcal B}$ is a Banach algebra with identity $I$. As usual, the set of bounded linear functionals on ${\mathcal B}$ is denoted by ${\mathcal B}^*$.
An element $\phi\in {\mathcal B}^*$ is said to be {\em central} if $\phi(XY)=\phi(YX)$ for all
$X,\ Y\in {\mathcal B}$. The set of central linear functionals on ${\mathcal B}$ shall be denoted by ${\mathcal B}_c^*$.
It is easy to see that ${\mathcal B}_c^*$ is a closed subspace of ${\mathcal B}^*$. A bounded linear functional $\phi$ on 
${\mathcal B}$ is said to be {\em multiplicative} if $\phi(XY)=\phi(X)\phi(Y)$. Clearly, multiplicative linear 
functionals are central. When ${\mathcal B}$ is commutative, the collection of multiplicative linear 
functionals is called the {\em maximal ideal space} of ${\mathcal B}$. 

If ${\mathcal B}$ is a $C^*$-algebra, then a positive central linear functional is called a {\em trace}. Of course, not every 
$C^*$-algebra possesses a trace.

Unless stated otherwise, $A=(A_0,\ A_1,\ ...,\ A_n)$ always stands for an $(n+1)$-tuple of general elements in 
${\mathcal B}$. A tuple $A$ is said to be commutative if $A_iA_j=A_jA_i$, $\forall 0\leq i,\ j\leq n$. 
In this paper, the ${\mathcal B}$-valued linear function $A(z)=z_0A_0+z_1A_1+\cdots z_nA_n$ is a primary associate of a tuple $A$. Without loss of generality, we assume 
the elements $A_0,\ A_1,\ ...,\ A_n$ are linearly independent, hence the range of $A(z)$ is an $n+1$ dimensional subspace of 
of ${\mathcal B}$, which we denote by $E_A$. A subalgebra ${\mathcal A}$ of ${\mathcal B}$ is said to be inversion-closed if
for every invertible element $a\in {\mathcal A}$, $a^{-1}$ is also in ${\mathcal A}$. For a tuple $A=(A_0,\ A_1,\ ...,\ A_n)$,
we let ${\mathcal B}_A$ denote the smallest inversion-closed Banach sub-algebra of ${\mathcal B}$ that 
contains $A_0,\ A_1,\ ...,\ A_n.$ Clearly, $A(z)$ is invertible in ${\mathcal B}$ if and only if it is invertible in 
${\mathcal B}_A$. Moreover, when $A$ is a commutative tuple, ${\mathcal B}_A$ is a commutative Banach algebra. In this case,
the maximal ideal space shall be denoted by $M_A$.\\
    
{\bf Definition.} For a tuple $A$, we let 
\[P(A)=\{z\in {\cc}^{n+1}:\ A(z)\ \text{is not invertible}.\}\]
The {\em projective spectrum} $p(A)$ of $A$ is $P(A)^T$, e.g.
\[p(A)=\{z=[z_0,\ z_1,\ ...,\ z_n]\in \pn:\ A(z)\ \text{is not invertible}.\}\]
For simplicity, we also refer to $P(A)$ as projective spectrum. The {\em projective resolvent sets} refer to
their complements $p^c(A)=\pn \setminus p(A)$ and $P^c(A)={\cc}^{n+1}\setminus P(A)$.\\

We let ${\mathcal B}^{-1}$ be the set of invertible
elements in ${\mathcal B}$. It is easy to see that the linear isomorphism $A(\cdot):\ \cc^{n+1}\longrightarrow E_A$ is a homeomorphism from
$P^c(A)$ to $E_A\cap {\mathcal B}^{-1}$, as well as from $p^c(A)$ to $(E_A\cap {\mathcal B}^{-1})/\cc^{\times}$.
Since $E_A\cap {\mathcal B}^{-1}$ is open in $E_A$, $P^c(A)$ is open, hence $p(A)$ is a 
compact subset of $\pn$. In some cases $p(A)$ can be equal to the entire space $\pn$, for example, when $A$ is a tuple of compact 
operators on an infinite dimensional Hilbert space. But one can always consider the slightly bigger tuple $(I,\ A_0,\ A_1,\ ...,\ A_n)$
which clearly has a more interesting projective spectrum. 
So without loss of generality, we assume throughout of the paper that $p(A)$ is a proper subset of $\pn$, or equivalently
$P(A)\neq \cc^{n+1}$.

First, we establish the nontriviality of $p(A)$. Idea of proof is from [Ya].
\begin{prop}
For any tuple $A$, $p(A)$ is a nontrivial compact subset of $\pn$.
\end{prop}
\begin{pf}
It only remains to show that $p(A)$ nontrivial, or equivalently,
$P(A)$ contains elements other than the origin $0$.

One first checks that on $P^c(A)$,
\begin{align*}
&A^{-1}(z_0,\ z_1,\ ...,\ z_n)-A^{-1}(z'_0,\ z_1,\ ...,\ z_n)\\
&=A^{-1}(z)\left(I-(z_0A_0+z_1A_1+\cdots z_nA_n)(z'_0A_0+z_1A_1+\cdots z_nA_n)^{-1}\right)\\
&=A^{-1}(z)\big(I-((z_0-z'_0)A_0+(z'_0A_0+z_1A_1+\cdots z_nA_n))(z'_0A_0+z_1A_1+\cdots z_nA_n)^{-1}\big)\\
&=-(z_0-z'_0)A^{-1}(z)A_0(z'_0A_0+z_1A_1+\cdots z_nA_n)^{-1}.
\end{align*}
This shows that $A^{-1}(z)$ is analytic in $z_0$, and likewise in all other variables. Moreover, the calculations show that
\begin{equation*}
\frac{\partial}{\partial z_j}A^{-1}(z)=-A^{-1}(z)A_jA^{-1}(z),\ \ \ 0\leq j\leq n. \tag{1.1}
\end{equation*} 
By Hartogs extension theorem, the origin $0$ cannot be an isolated singularity of $A^{-1}(z)$, and hence $p(A)$ is nontrivial.
\end{pf}

Proposition 1.1 is an interesting fact, since the elements $A_0,\ A_1,\ ...,\ A_n$ may have nothing to do with each other.
Let us look at a few examples.\\

{\bf Example 1}. When ${\mathcal B}$ is the matrix algebra $M_k(\cc)$, $A=(A_0,\ A_1,\ ...,\ A_n)$ is a tuple of $k\times k$ matrices.
Then $A(z)$ is invertible if and only if $detA(z)\neq 0$. Since $detA(z)$ is homogenous of degree $k$,
\[p(A)=\{z=[z_0,\ ,z_1,\ ...,\ z_n]\in \pn:\ detA(z)=0\}\]
is a projective hypersurface of degree $k$. And $p^c(A)$ in this case is a {\em hypersurface complement}.\\

{\bx}

{\bf Example 2}. Let $A_0$ be any element in ${\mathcal B}$, and $A_1=-I$. Then for the tuple $A=(A_0,\ A_1)$,
$A(z)=z_0A_0-z_1I$. Clearly, if $[z_0,\ z_1]$ is in $p(A)$
then $z_0\neq 0$, and $p(A)$ under the affine coordinate $z_1/z_0$ is the classical spectrum $\sigma(A_0)$. So Proposition 1.1 in fact implies the 
nonemptiness of the classical spectrum.\\

{\bx}

{\bf Example 3}. Now consider $L^2({\mathbb T},\ m)$, where $m$ is the normalized Lebesgue measure on the unit circle $\td$. 
$\{w^n:\ \ n\in \z,\ |w|=1\}$ is an orthonormal basis for $L^2({\mathbb T},\ m)$. Let $\theta$ be an irrational number and 
set $\lb=exp(2\pi \sqrt{-1}\theta)$. Consider the two unitaries defined by 
\[A_0f(w)=wf(w),\ \ \ A_1f(w)=f(\lb w),\ \ f\in L^2({\mathbb T},\ m),\]
and let ${\mathcal B}$ be the $C^*$-algebra generated by $A_0$ and $A_1$. 
Clearly, $A(z)$ is invertible if and only if $z_0A_0A^*_1+z_1I$ is invertible. So by Example 2,
$p(A)=-\sigma(A_0A_1^*)$. One checks that $A_1^*f(w)=f({\bar \lb}w)$, hence
\[A_0A_1^*w^n=w({\bar \lb}w)^n={\bar \lb}^nw^{n+1}.\]
So $A_0A_1^*$ is a unitary bilateral weighted shift. Since $A_1A_0=\lb A_0A_1$,
\[A_1^*A_0A_1^*A_1=A_1^*A_0={\bar \lb}A_0A_1^*,\]
hence $\sigma(A_0A_1^*)$ is invariant under multiplication by ${\bar \lb}$ which implies $\sigma(A_0A_1^*)=\td$. Therefore
\[p(A)=\td,\ \ \  \text{and}\ \ \  P(A)=\{(z_0,\ z_1)\in \cc^2:\ |z_0|=|z_1|\}.\]
In this case, $P^c(A)$ consists of two connected components:
\[\Omega_0=\{ |z_0|>|z_1|\},\ \ \text{and}\ \ \Omega_1=\{ |z_0|<|z_1|\}.\] 

The $C^*$ algebra generated by $A_0$ and $A_1$ is the {\em irrational rotation algebra} often denoted by 
${\mathcal A}_{\theta}$. We will come back to this algebra in Section 3.\\

{\bx}

\section{Projective spectrum and projective hypersurface}

As we have seen in Example 1, for ${\mathcal B}=M_k(\cc)$, a projective spectrum is a projective hypersurface. Naturally, things could become very 
different for other Banach algebras. For instance, in Example 3 the projective resolvent sets are disjoint unions of two connected 
components, while a hypersurface complement is always connected. Nevertheless, as it turns out, projective spectrum resembles hypersurface 
in many other ways.

In $\pn$ for $n\geq 2$, a {\em line} is the quotient of a two dimensioal subspace of $\cc^{n+1}$ (removing the origin) over $\cc^{\times}$.
By virtue of the Fundamental Theorem of Algebra, a hypersurface in $\pn$ intersects with every line. This fact holds true 
for $p(A)$.

\begin{corr}
For $n\geq 2$, every line in $\pn$ intersects with $p(A)$.
\end{corr}
\begin{pf}
It is equivalent to show that every two dimensional subspace in $\cc^{n+1}$ intersects with $P(A)$ nontrivially. In fact,
for any two linearly independent vectors $\lb=(\lb_0,\ \lb_1,\ ...,\ \lb_n)$ 
and $\eta=(\eta_0,\ \eta_1,\ ...,\ \eta_n)$ in $\cc^{n+1}$, let
\[A'=\sum_{j=0}^n \lb_jA_j,\ \ A''=\sum_{j=0}^n\eta_jA_j.\]
By Proposition 1.1 for the case $A=(A',\ A'')$, there exists scalars $a$ and $b$, not both zero, such that
$aA'+bA''$ is not invertible, and hence $a\lb+b\eta \in P(A)$.
\end{pf}

In the case $A$ is a commutative tuple, the projective spectrum can be explicitly calculated. 

\begin{prop}
If $A$ is a commutative tuple, then $P(A)$ is a union of hyperplanes.
\end{prop}
\begin{pf}
As remarked in Section 1, ${\mathcal B}_A$ in this case is commutative, and $P(A)$ is unchanged when considered in 
${\mathcal B}_A$. Let $M_A$ denote the maximal ideal space of ${\mathcal B}_A$. Then by Gelfand theorem $A(z)$ is not 
invertible in ${\mathcal B}_A$ if and only if there exists $\phi\in M_A$ such that 
\begin{equation*}
\phi(A(z))=\sum_{j=0}^{n}z_j\phi(A_j)=0.\tag{2.1}
\end{equation*}
For simplicity, we let $H_\phi=\{z\in \cc^{n+1}:\ \sum_{j=0}^{n}z_j\phi(A_j)=0\}.$ 
If $\phi$ is such that $\phi(A_j)=0$ for all $j$, then $P(A)=\cc^{n+1}$, which is clearly a union 
(uncountable) of hyperplanes. Otherwise $H_\phi$ is a hyperplane, and one sees that
\begin{equation*}
P(A)=\cup_{\phi\in M_A}H_{\phi}.\tag{2.2}
\end{equation*}
\end{pf}

In the case when $M_A$ is a finite set, for instance when ${\mathcal B}=M_k(\cc)$, $P(A)$ is a union of a finite number 
of hyperplanes, e.g., $P(A)$ is a {\em central hyperplane arrangement}. In this case the topology of $P^c(A)$ is a primary topic in 
Hyperplane Arrangement (cf. Orlik and Terao [OT]). It is worth mentioning that every central hyperplane arrangement can be represented as the projective spectrum of 
a tuple of diagonal matrices. For example, for the {\em braid arrangement} in $\cc^3$ defined by equation $(z_0-z_1)(z_1-z_2)(z_2-z_0)=0$,
one can let 
\[
A_0=\begin{pmatrix}
1 & 0 & 0\\
0 & -1 & 0\\
0 & 0 & 0 \\
\end{pmatrix},\ \ \ 
A_1=\begin{pmatrix}
-1 & 0 & 0\\
0 & 0 & 0\\
0 & 0 & 1 \\
\end{pmatrix},\ \ \ 
A_2=\begin{pmatrix}
0 & 0 & 0\\
0 & 1 & 0\\
0 & 0 & -1 \\
\end{pmatrix},\]
and verifies easily that $P(A)=\{z\in \cc^3:\ (z_0-z_1)(z_1-z_2)(z_2-z_0)=0\}$.\\

{\bf Example 4}. Now consider the disk algebra ${\mathcal B}=A({\bf D})$ and let $A_j=w^j,\ 0\leq j\leq n$. Then
$A$ is a commutative tuple, and 
\[A(z)=\sum_{j=0}^n z_jw^j\]
is a degree $n$ polynomial in $w$ (when $z_n\neq 0$). In this case, $A(z)$ is invertible in ${\mathcal B}$ if and only if it has no zero in the 
closed unit disk ${\overline {\bf D}}$. Here, the maximal ideal space of ${\mathcal B}$ is equal to ${\overline {\bf D}}$ (cf. Douglas [Do]), and a point 
$w\in {\overline {\bf D}} $ acts on ${\mathcal B}$ by point evaluation 
\[\phi_w(f):=f(w),\ \ f\in {\mathcal B}.\]
Then by Proposition 2.2.
\[P(A)=\bigcup_{|w|\leq 1}H_{\phi_{w}}.\]

{\bx}

When ${\mathcal B}=M_k(\cc)$, $P(A)$ is a hypersurface in $\cc^{n+1}$ defined by the polynomial equation $detA(z)=0$. So its complement is clearly a domain of holomorphy, because
$1/detA(z)$ is holomorphic on $P^c(A)$ and cannot be extended analytically to a neighborhood of any point in $P(A)$. On another account,
one direct consequence of Proposition 2.2 is that when $A$ is commutative, each path connected component of $P^c(A)$ is a domain of 
holomorphy. To see this, we let $U$ be a connected component of $P^c(A)$, and $\lb$ be any point in $\partial U$.
Since $P(A)$ is a union of hyperplanes, $\lb$ is in one of these hyperplanes, say $\{\sum_{j=0}^na_jz_j=0\}$. So
$(\sum_{j=0}^na_jz_j)^{-1}$ is holomorphic on $U$ and does not have a holomorphic extension to any neighborhood of $\lb$.

On a general projective resolvent set $P^c(A)$, $A^{-1}(z)$ is holomorphic and cannot be extended to any greater region. 
 So it is natural to ask whether $P^c(A)$
is necessarily a domain of holomorphy, or a disjoint union of domains of holomorphy when it is not path connected. Of course,
the difference here is that $A^{-1}(z)$ is a ${\mathcal B}$-valued function. 

Here, we show that for some interesting types of Banach algebra ${\mathcal B}$, the answer is positive.

\begin{thm}
If ${\mathcal B}$ is reflexive as a Banach space, i.e. ${\mathcal B}={\mathcal B}^{**}$, then every connected component of $P^c(A)$ 
is a domain of holomorphy.
\end{thm}
\begin{pf}
We let $U$ be a connected component of $P^c(A)$, and $\lb$ be a point in $\partial U$. We will show by contracdiction that there exists
a $\phi\in {\mathcal B}^*$ such that $\phi(A^{-1}(z))$ does not extend holomorphicly to any neighborhood of $\lb$.

Suppose on the contrary for every $\phi\in {\mathcal B}^*$, $\phi(A^{-1}(z))$ extends holomorphicly to a neighborhood of $\lb$. Then one 
observes that the function 
\[F(\phi,\ z):=\phi(A^{-1}(z)),\ \ \phi\in {\mathcal B}^{*},\ z\in U,\]
is a bounded linear functional on ${\mathcal B}^*$ for every fixed $z$, and has, for every fixed $\phi$, a holomorphic continuation to 
a neighborhood of $\lb$.
Let $z^m,\ m\geq 0,$ be a sequence in $U$ that converges to $\lb$, and consider the sequence $F_m(\phi):=F(\phi,\ z^m)$.
To take care of the case that $\partial U$ may intersect itself at $\lb$, we assume that for every open neighborhood $V$ 
of $\lb$, $z^m$ stay in the same connected component of $V\cap U$ when $m$ is sufficiently large. 

Then $F_m\in {\mathcal B}^{**},\ \forall m,$ and for every fixed $\phi$
\[sup\{|F_m(\phi)|:\ m\geq 0\}<\infty.\]
The Uniform Boundedness Principle then implies that the limit
\[F_{\infty}(\phi):=\lim_{m\to \infty}F_m(\phi),\ \ \phi\in {\mathcal B}^*\]
is in ${\mathcal B}^{**}$. Since ${\mathcal B}={\mathcal B}^{**}$, there exists a $B\in {\mathcal B}$ such that
\begin{equation*}
F_{\infty}(\phi)=\phi(B),\ \ \ \forall \phi\in {\mathcal B}^*.\tag{2.3}
\end{equation*}

Moreover, for a fixed $C\in {\mathcal B}$ and any $\phi\in {\mathcal B}^*$, the functional $\phi_C$ defined by
\[\phi_C(X):=\phi(XC),\ \ X\in {\mathcal B}\]
is clearly in ${\mathcal B}^*$, so it follows from (2.3)
that 
\begin{align*}
\lim_{m\to \infty}\phi(A^{-1}(z^m)C)&=\lim_{m\to \infty}F(\phi_c,\ z^m)\\
&=\lim_{m\to \infty}F_m(\phi_c)\\
&=F_{\infty}(\phi_c)\\
&=\phi(BC), \ \ \ \forall \phi\in {\mathcal B}^*.
\end{align*}
Letting $C=A(\lb)$, we have
\begin{align*}
\phi(BA(\lb))&=\lim_{m\to \infty}\phi(A^{-1}(z^m)A(\lb))\\
&=\lim_{m\to \infty}\phi\left((A^{-1}(z^m)(A(z^m)+A(\lb)-A(z^m)))\right)\\
&=\phi(I)+\sum_{j=0}^n\lim_{m\to \infty}(\lb_j-z^m_j)\phi(A^{-1}(z^m)A_j)\\
&=\phi(I)+\sum_{j=0}^n\lim_{m\to \infty}(\lb_j-z^m_j)\phi_{A_j}(A^{-1}(z^m)).
\end{align*} 
Since $\phi_{A_j}(A^{-1}(z))$ extends analytically to a neighborhood of $\lb$, $\phi_{A_j}(A^{-1}(z^m))$ is bounded, and it follows that
\[\phi(BA(\lb))=\phi(I).\ \ \ \forall \phi\in {\mathcal B}^*. \]
 
Similarly, we can also show $\phi(A(\lb)B)=\phi(I),\ \forall \phi\in {\mathcal B}^*$. These imply that
$BA(\lb)=A(\lb)B=I$, e.g, $A(\lb)$ is invertible, which is a contradiction.
\end{pf}

The proof of Theorem 2.3 can be modified to work for other Banach algebras. For example, if ${\mathcal H}$ is a reflexive Banach space,
and ${\mathcal B}$ is a Banach sub-algebra of $B({\mathcal H})$---the set of bounded linear operators on ${\mathcal H}$, then for every
$x\in {\mathcal H}$ and $f\in {\mathcal H}^*$ , \[\phi_{x,f}(C)=f(Cx),\ \ C\in {\mathcal B}\] defines a bounded linear functional on 
${\mathcal B}$. If we let $F_m(x,\ f):=f(A^{-1}(z^m)x)$ and apply the Uniform Boundedness Principle, then
\[F_{\infty}(x,\ f):=\lim_{m\to \infty}F_m(x,\ f),\ \ \ \  x\in {\mathcal H}, f\in {\mathcal H}^*\]
is a bounded bilinear form on ${\mathcal H}\times {\mathcal H}^*$. In particular, if we fix $x$ then $F_{\infty}(x,\ \cdot)$
is in ${\mathcal H}^{**}$. Now since ${\mathcal H}$ is reflexive, there is a unique $B(x)\in {\mathcal H}$ such that
\[F_{\infty}(x,\ f)=f(B(x)),\ \forall x\in {\mathcal H},\ f\in {\mathcal H}^*,\]
and it is not hard to see that $B$ is a linear and bounded.
 
Similar to the ending part of the proof of Theorem 2.3, we have
\[f(BA(\lb)x)=f(x)=f(A(\lb)Bx),\ \ \forall x\in {\mathcal H},\ f\in {\mathcal H}^*,\]
which means $A(\lb)$ is invertible with inverse $B$. However, in general this $B$ may not be in ${\mathcal B}$.

But things can be pulled together in the case when ${\mathcal B}$ is a $C^*$-algebra. In this case, ${\mathcal B}$ can be 
identified (up to a isometrically $*$-isomorphism) with a $C^*$-subalgebra of $B({\mathcal H})$ (cf. Davidson [Da]), where ${\mathcal H}$ is a Hilbert 
space. And an element in ${\mathcal B}$ is invertible if and only if it is invertible in $B({\mathcal H})$ (cf. Douglas [Do]). We therefore have the following

\begin{thm}
If ${\mathcal B}$ is a unital $C^*$-algebra, then every connected component of $P^c(A)$ is a domain of holomorphy.\\
\end{thm}  

{\bf Question A}. Is the statement in Theorem 2.4 true for any unital Banach algebra?\\

\section{${\mathcal B}$-valued 1-form $A^{-1}(z)dA(z)$ and the de Rham cohomology space $H^1_d(P^c(A),\ \cc)$}

If $S$ is a hypersurface defined by $\{q(z)=0\}$, where $q$ is an irreducible homogenous polynomial of degree $k>0$,
then the complements $S^c=\cc^{n+1}\setminus S$ and $(S^T)^c=\pn\setminus S^T$ are both connected. Moreover,
the singular homology group $H_1(S^c,\ \z)=\z$, and $H_1((S^T)^c,\ \z)=\z/k\z$ (cf. Dimca [Di], Ch4), which indicates that
neither $S^c$ nor $(S^T)^c$ is simply connected. But, as indicated in Example 3, projective resolvent sets $P^c(A)$ and 
$p^c(A)$ may not be connected. Furthermore, connected components of $p^c(A)$ may also be simply connected.
However, connected components of $P^c(A)$ behave somewhat differently. In this section we will have a peek on the 
de Rham cohomology space $H^1_d(P^c(A),\ \cc)$. The Maurer-Cartan type ${\mathcal B}$-valued 1-form 
$\omega_A(z):=A^{-1}(z)dA(z)$ and central linear functionals on ${\mathcal B}$ are important tools in our study.
Here $d=\sum_{j=0}^n\frac{\partial}{\partial z_j}dz_j.$

First of all, from a operator-theoretic point of view, $\omega_A(z)$ is a faithful associate of the tuple $A$ because it determines
$A$ up to a certain equivalence. Let $A$ and $B$ be two tuples with the same projective spectrum $S\subset \cc^{n+1}$. Here, the two 1-forms $\omega_A(z)$ and 
$\omega_B(z)$ are said to be {\em similar} if there is an invertible element $V\in {\mathcal B}$ such that
\[V^{-1}\omega_A(z)V=\omega_B(z),\ \ \ \forall z\in S^c.\]

\begin{prop}
Let $A$ and $B$ be two tuples with the same projective spectrum $S$. Then the two 1-forms $\omega_A(z)$ and 
$\omega_B(z)$ are similar if and only if there are invertible element $U,\ V\in {\mathcal B}$ such that
$UA_jV=B_j$ for every $0\leq j\leq n$.
\end{prop}
\begin{pf}
For the sufficiency, one easily checks that $UA_jV=B_j$ for each $0\leq j\leq n$ implies $UA(z)V=B(z)$. Then on the projective resolvent set
$S^c$,
\begin{align*}
B^{-1}(z)dB(z)&=V^{-1}A^{-1}(z)U^{-1}dUA(z)V\\
&=V^{-1}A^{-1}(z)U^{-1}UdA(z)V\\
&=V^{-1}A^{-1}(z)dA(z)V.
\end{align*}

For the necessity, one checks that $V^{-1}\omega_A(z)V=\omega_B(z)$ implies
\[\sum_{j=0}^nV^{-1}A^{-1}(z)A_jVdz_j=\sum_{j=0}^nB^{-1}(z)B_jdz_j,\]
and hence $V^{-1}A^{-1}(z)A_jV=B^{-1}(z)B_j$ for each $j$, or equivalently,
\begin{equation*}
A_jV=A(z)VB^{-1}(z)B_j,\ \ \ \forall z\in S^c. \tag{3.1}
\end{equation*}

So for any fixed $w\in S^c$, one has
\[\sum_{j=0}^nw_jA_jV=A(z)VB^{-1}(z)\sum_{j=0}^nw_jB_j,\]
which implies that
\[A(w)VB^{-1}(w)=A(z)VB^{-1}(z),\ \ \ \forall z\in S^c.\]
So $A(z)VB^{-1}(z)$ is an invertible constant, for which we denote by $U^{-1}$. Then by (3.1), $UA_jV=B_j$ for every $0\leq j\leq n$.
\end{pf}  

One observes that for a $\phi\in {\mathcal B}^*$, $\phi(\omega_A(z))=\sum_{j=0}^n\phi(A^{-1}(z)A_j)dz_j$ is a holomorphic 1-form
on $P^c(A)$. 
\begin{thm}
Consider a bounded linear functional $\phi$ on ${\mathcal B}$.

(a) If $\phi$ is central, then $\phi(\omega_A(z))$ is a closed 1-form on $P^c(A)$. 

(b) If $\phi(I)\neq 0$, then there is no holomorphic function $f(z)$ on $P^c(A)$ such that
$df(z)=\phi(\omega_A(z))$.
\end{thm}
\begin{pf}
(a) First of all, Maurer-Cartan type form $\omega$ has the property $d\omega=-\omega\wedge \omega$. This fact for $\omega_A(z)$  
also follows easily from (1.1). By (1.1), for each $j$,
\[\frac{\partial}{\partial z_j}A^{-1}(z)=-A^{-1}(z)A_jA^{-1}(z),\]
hence
\begin{align*}
d\omega_A(z)&=\sum_{j=0}^ndA^{-1}(z)A_jdz_j\\
&=\sum_{i,j=1}^n\frac{\partial A^{-1}(z)}{\partial z_i}A_jdz_i\wedge dz_j\\
&=\sum_{i,j=0}^n-A^{-1}(z)A_iA^{-1}(z)A_jdz_i\wedge dz_j\\
&=\sum_{i<j}-(A^{-1}(z)A_iA^{-1}(z)A_j -A^{-1}(z)A_jA^{-1}(z)A_i)dz_i\wedge dz_j \tag{3.2}\\
&=-\omega_A(z)\wedge \omega_A(z). \tag{3.3}
\end{align*}
If $\phi$ is central, then by (3.2)
\begin{align*}
d\phi(\omega_A(z))&=\phi(d\omega_A(z))\\
&=\sum_{i<j}-\phi(A^{-1}(z)A_iA^{-1}(z)A_j -A^{-1}(z)A_jA^{-1}(z)A_i)dz_i\wedge dz_j\\
&=0,
\end{align*}
hence $\phi(\omega_A(z))$ is closed.

(b) If there exists an $f$ holomorphic on $P^c(A)$ such that
$df(z)=\phi(\omega_A(z))$, then 
\[\frac{\partial f}{\partial z_j}(z)=\phi(A^{-1}(z)A_j),\ \ \forall j.\]
So for any $t\in \cc^{\times }$,
\begin{equation*}
\frac{\partial f}{\partial z_j}(tz)=\phi(A^{-1}(tz)A_j)=t^{-1}\phi(A^{-1}(z)A_j),\ \ \forall j.\tag{3.4}
\end{equation*}
It follows that 
\begin{align*}
d(f(tz))&=\sum_{j}t\frac{\partial f}{\partial z_j}(tz)dz_j\\
&=tt^{-1}\sum_{j}\phi(A^{-1}(z)A_j)dz_j\\
&=\phi(\omega_A(z))\\
&=df(z),
\end{align*}
hence $f(tz)-f(z)$ is a constant depending on $t$, say $c(t)$. To figure out $c(t)$, one computes using (3.4) that
\begin{align*}
c'(t)&=\frac{df(tz)}{dt}\\
&=\sum_{j}z_j\frac{\partial f}{\partial z_j}(tz)\\
&=t^{-1}\sum_{j}z_j\phi(A^{-1}(z)A_j)\\
&=t^{-1}\phi(I).
\end{align*}
Since $c(1)=0$, $c(t)=\phi(I)logt$, hence
\begin{equation*}
f(tz)-f(z)=\phi(I)logt.\tag{3.5}
\end{equation*}
Since $f$ is holomorphic on $P^c(A)$, and $tz\in P^c(A)$ for $t\in \cc^{\times}$, $f(tz)$ is holomorphic in $t$, and hence
$\phi(I)logt$ is holomorphic on $\cc^{\times}$, which is possible only if $\phi(I)=0$.
\end{pf}

If $P^c(A)$ is not path connected, Theorem 3.2 can be stated for every connected component of $P^c(A)$.
For a domain $U\in \cc^{n+1}$ of holomorphy ( or equivalently, a {\em  Stein domain}), its de Rham cohomology $H^*_d(U,\ \cc)$ 
can be calculated by holomorphic forms (cf. Range [Ra]). To be precise, if $\Omega^r(U)$ is the collection of 
holomorphic $r$-forms on $U$, then
\[H^r_d(U,\ \cc)\simeq \{f\in \Omega^r(U):\ df=0\}/d\Omega^{r-1}(U),\ \ \ for\ r\geq 0.\]   
This observation, combined with Theorem 2.4 and Theorem 3.2, leads to the following
\begin{thm}
Let ${\mathcal B}$ be a $C^*$-algebra with a trace $\phi$, and $U$ be a connected component of $P^c(A)$.
Then $\phi(\omega_A(z))|_U$ is a nontrivial element in the de Rham 
cohomology space $H^1_d(U,\ \cc)$. In particlar, $U$ is not simply connected.
\end{thm}

\begin{pf}
By Theorem 2.4, $U$ is a domain of holomorphy. Hence 
the de Rham cohomology space $H^1_d(U,\ \cc)$ can be calculated by holomorphic forms.
Now since $\phi$ is a trace, $\phi $ is central with $\phi(I)>0$. Theorem 3.2 then concludes that $\phi(\omega_A(z))|_U$ is a nontrivial element in $H^1_d(U,\ \cc)$.\\  
\end{pf}

{\bf Example 5}. Let $A$ be a tuple of $k\times k$ matrices, and let $Tr$ be the ordinary trace on square matrices.
$P^c(A)$ is a hypersurface complement in this case.
It is a classical fact that for a one variable square matrix-valued differentiable function $M(t)$
\[Tr(M^{-1}(t)\frac{d}{dt}M(t))=\frac{d}{dt}\text{log(det}M(t)),\]
hence we have
\[Tr(\omega_A(z))=d\text{logdet}A(z),\ \ \ z\in P^c(A).\]
One sees that $\text{logdet}A(z)$ is not a global holomorphic function on $P^c(A)$, hence $Tr(\omega_A(z))$ is 
closed but not exact.\\ 

{\bx} 

In order to have more traces, one can let ${\mathcal B}={\mathcal B}_A$.\\
 
{\bf Example 6}.  Consider a tuple $A$ of $3\times 3$ matrices, where 
\[
A_0=\begin{pmatrix}
1 & 0 & 0\\
0 & 0 & -1\\
0 & 1 & 0 \\
\end{pmatrix},\ \ \ 
A_1=\begin{pmatrix}
1 & 0 & 0\\
0 & \sqrt{-1} & 0\\
0 & 0 & -\sqrt{-1} \\
\end{pmatrix},\ \ \ 
A_2=\begin{pmatrix}
1 & 0 & 0\\
0 & 1 & 0\\
0 & 0 & 1 \\
\end{pmatrix}.\]
Then 
\[A(z)=
\begin{pmatrix}
z_0+z_1+z_2 & 0 & 0\\
0 & \sqrt{-1}z_1+z_2 & -z_0\\
0 & z_0 & -\sqrt{-1}z_1+z_2 \\
\end{pmatrix},\]
and $detA(z)=(z_0+z_1+z_2)(z_0^2+z_1^2+z_2^2)$.
Hence \[P(A)=\{(z_0,\ z_1,\ z_2)\in {\cc}^3:\ (z_0+z_1+z_2)(z_0^2+z_1^2+z_2^2)=0\}.\] 

In this case, one verifies that ${\mathcal B}_A=\cc\oplus M_{2\times 2}(\cc)$.
Let $\phi_1$ and $\phi_2$ be the linear functionals on $\cc\oplus M_{2\times 2}(\cc)$ defined by
\[\phi_1 
\begin{pmatrix}
a & 0\\
0 & T\\
\end{pmatrix}=a,\ \ \ 
\phi_2 
\begin{pmatrix}
a & 0\\
0 & T\\
\end{pmatrix}=Tr(T),\] 
where $T\in M_{2\times 2}(\cc)$. Then $\phi_1$ and $\phi_2$ are both traces, hence by Example 5
\[\phi_1(\omega_A(z))=\frac{dz_0+dz_1+dz_2}{z_0+z_1+z_2},\ \ \ 
\phi_2(\omega_A(z))=\frac{2z_0dz_0+2z_1dz_1+2z_2dz_2}{z_0^2+z_1^2+z_2^2}\]
are nontrivial 1-forms in $H^1_d(P^c(A),\ \cc)$. \\

{\bx}

{\bf Example 7}. Now we continue with Example 3. We have remarked that is this case ${\mathcal B}$ is the irrational
rotation algebra ${\mathcal A}_{\theta}$. It is not hard to see that 
monomials $A_0^kA_1^l,\ \ k,\ l\in {\mathbb Z},$ span a dense subspace of ${\mathcal A}_{\theta}$. 
${\mathcal A}_{\theta}$ possesses a unique faithful unital trace $\phi$ defined by $\phi(I)=1$ and 
$\phi(A_0^kA_1^l)=0$ for $k$ and $l$ not both zero. 

We now compute $\phi(\omega_A(z))$, where $A=(A_0,\ A_1)$. We showed in Example 3 that in this case $P^c(A)$ has two 
connected components:
\[\Omega_0=\{ |z_0|>|z_1|\},\ \ \text{and}\ \ \Omega_1=\{ |z_0|<|z_1|\}.\] 
On $\Omega_0$,
\begin{align*}
A^{-1}(z)dA(z)&=(z_0A_0+z_1A_1)^{-1}(A_0dz_0+A_1dz_1)\\
&=z_0^{-1}(A_0+\frac{z_1}{z_0}A_1)^{-1}A_0(dz_0+A_0^{-1}A_1dz_1)\\
&=z_0^{-1}(I+\frac{z_1}{z_0}A_0^{-1}A_1)^{-1}(dz_0+A_0^{-1}A_1dz_1)\\
&=\left(I-\frac{z_1}{z_0}A_0^{-1}A_1+(\frac{z_1}{z_0}A_0^{-1}A_1)^2-\cdots \right)(\frac{dz_0}{z_0}+A_0^{-1}A_1\frac{dz_1}{z_0}).
\end{align*}
Hence $\phi(\omega_A(z))=\frac{dz_0}{z_0}$. Likewise, on $\Omega_1$, $\phi(\omega_A(z))=\frac{dz_1}{z_1}$.

As a matter of fact, in this case it is not hard to compute directly that
\[H^1_d(U_0,\ \cc)=\cc \frac{dz_0}{z_0},\ \ \ H^1_d(U_1,\ \cc)=\cc \frac{dz_1}{z_1}.\]\\

{\bx}
  
We conclude this section with a remark on the value of $\phi(I)$. The case $\phi(I)\neq 0$ is important for Theorem 3.2(b), and
as we will see, the case $\phi(I)=0$ is also meaningful.

For $\phi\in {\mathcal B}^*$,
one easily sees that $\phi(\omega_A(z))$ is homogenous of degree $0$. So it makes one wonder for what type of $\phi$,
$\phi(\omega_A(z))$ defines a 1-form on the projective resolvent set $p^c(A)\subset \pn$. 
For a holomorphic 1-form $\eta=\sum_{j=0}^{n}f_j(z)dz_j$, its {\em contraction} with the Euler
field $\theta=\sum_{j=0}^{n}z_j\frac{\partial}{\partial z_j}$ is
\[\Delta(\eta):=\sum_{j=0}^{n}z_jf_j(z).\] 
By Griffiths [Gr], a locally holomorphic 1-form $\eta=\sum_{j=0}^{n}f_j(z)dz_j$ on $\cc^{n+1}$ comes from 
a 1-form on $\pn$ if and only if it is homogenous of degree $0$ and the contraction $\Delta (\eta)=0$.   

It is easy to see that the Maurer-Cartan form $\omega_A(z)$ is homogenous of degree $0$.
But it itself is not a 1-form on $p^c(A)$. In fact, one checks easily that
\begin{align*}
\Delta(\omega_A(z))&=\Delta(\sum_{j=0}^nA^{-1}(z)A_jdz_j)\\
&=\sum_{j=0}^nA^{-1}(z)z_jA_j\\
&=I.
\end{align*}

So if $\phi$ is a linear functional on ${\mathcal B}$ such that $\phi(I)=0$, then 
\begin{align*}
\Delta(\phi(\omega_A(z)))&=\phi(\Delta(\omega_A(z)))\\
&=\phi(I)\\
&=0,
\end{align*}
hence $\phi(\omega_A(z))$ defines a global holomorphic 1-form on $p^c(A)$. Furthermore, using arguments similar to that in
the proof of Theorem 3.2(a), one can easily check that when, in addition, $\phi$ is central, $\phi(\omega_A(z))$ is also closed
on $p^c(A)$. 

\section{the case when $A$ is commutative}
  
When $A$ is a commutative tuple, ${\mathcal B}_A$ is an commutative Banach sub-algebra of ${\mathcal B}$. 
One observes that in this case, $({\mathcal B}_A)^*_c={\mathcal B}_A^*$, 
and every $\phi\in M_A$ has the property $\phi(I)=1$. Then, as remarked after Example 4, 
every connected component of $P^c(A)$ is a domain of holomorphy, and the next corollary is another consequence of Theorem 3.2. 

\begin{corr}
If $A$ is a commutative tuple, then for every $\phi\in M_A$, $\phi(\omega_A(z))$ is a nontrivial element in $H^1_d(P^c(A),\ \cc)$
\end{corr}

As stated in Proposition 2.2 that in this case $P(A)$ is a (possibly uncountable) union of hyperplanes. This section recalls
a theorem in Hyperplane Arrangements by Arnold and Briskorn, and discusses its possible analogue in the setting here.  

First, if $\phi$ is a multiplicative linear functional on ${\mathcal B}_A$, then
 \begin{align*}
\phi(\omega_A(z))&=\phi(A^{-1}(z))d\phi(A(z))\\
&=\frac{d\phi(A(z))}{\phi(A(z))}\\
&=\frac{d\sum_{j=0}^n z_j\phi(A_j)}{\sum_{j=0}^n z_j\phi(A_j)}\tag{3.6}
\end{align*}
Here one recalls that $\sum_{j=0}^n z_j\phi(A_j)$ is the defining function for the hyperplane $H_{\phi}$.

In the case $A$ is a tuple of generic commutative $k\times k$ matrices, the maximum ideal space $M_A$ consists of $k$ elements, say, $\phi_1,\ ...,\ \phi_k$, and hence 
\[P(A)=\cup_{j=0}^k H_{\phi_j}\]
is a central arrangement. By a well-known result conjectured by Arnold [Ar] and proved by Briskorn [Br], 
the cohomology algebra $H^*_d(P^c(A),\ \cc)$ is generated by $1$ and the 1-forms 
\[\frac{d\sum_{j=0}^n z_j\phi(A_j)}{\sum_{j=0}^n z_j\phi(A_j)},\ \ \ \phi\in M_A.\]

Here we make two observations.\\

{\bf 1}. First, let $\wedge ({\mathcal B}_A^*)$ be the exterior algebra on ${\mathcal B}_A^*$.
By Theorem 3.2, $\omega_A(z)$ induces a homomorphism 
$\omega_A(z)^*:\wedge ({\mathcal B}_A^*)\longrightarrow H^*_d(P^c(A),\ \cc)$ defined by $\omega_A(z)^*(1)=1$, and
\[\omega_A(z)^*(\phi_1\wedge \phi_2\wedge \cdots \wedge \phi_q)=
\phi_1(\omega_A(z))\wedge \phi_2(\omega_A(z))\wedge \cdots \wedge \phi_q(\omega_A(z)),\ \ q\geq 1.\]

Corollary 4.1 and the remarks above lead to the following simple consequence of the Arnold and Briskorn's theorem.
\begin{corr}
For ${\mathcal B}=M_k(\cc)$ and any tuple $A$ of commutative matrices in $M_k(\cc)$, 
$\omega_A(z)^*:\ \wedge ({\mathcal B}_A^*) \longrightarrow H^*(P^c(A),\ \cc)$ is surjective.
\end{corr}
Corollary 4.1 indicates that for ${\mathcal B}=M_k(\cc)$, topological information of $P^c(A)$ is all encoded in $\omega_A(z)$.
It will be interesting to see if there are other Banach algebras with this property.\\
    
{\bf 2}. For a general commutative tuple $A$, $P^c(A)$ may not be connected, hence $\omega_A(z)^*$ may not be surjective. But it makes sense to ask 
whether there is a similar fact for every connected component $U$ of $P^c(A)$, e.g., whether the map
$\omega_A(z)^*:\ \wedge ({\mathcal B}_A^*) \longrightarrow H^*(U,\ \cc)$ is surjective. But since in general both 
${\mathcal B}_A^*$ and $H^*(U,\ \cc)$ can be infinite dimensional, the nature of this problem is somewhat hard to see.

A related particular case is when the maximal ideal space $M_A$ is path connected. For any $\phi_0,\ \phi_1\in M_A$, let $\phi_t,\ t\in [0,\ 1]$ be a 
continuous path in $M_A$. Taking any fixed cycle $\gamma \in H_1(U, \z)$ and using (4.1), we have
\[\frac{1}{2\pi \sqrt{-1}}\int_{\gamma}\phi_t(\omega_A(z))=\frac{1}{2\pi \sqrt{-1}}\int_{\gamma}dlog(\sum_{j=0}^n z_j\phi_t(A_j)),\]
which is an integer-valued continuous function in $t$, hence is a constant. This means that 
$\omega_A(z)^*$ restricted to $M_A$ is a constant in $H^1_d(U,\ \cc)$.\\ 
 
{\bf Example 8.} Now we take another look at Example 4. Note that in this case ${\mathcal B}_A=A({\bf D})$. Then,
\[\omega_A(z)=\frac{\sum_{j=0}^nw^jdz_j}{\sum_{j=0}^nz_jw^j}.\]
Now $P^c(A)$ is connected (as we will see in a minute), and $M_A$ can be identified with $\overline{D}$ through evaluation(cf. Douglas [Do]), 
so as remarked above $\omega_A(z)^*$ restricted to $M_A$ 
is a constant in $H^1_d(P^c(A),\ \cc)$. 
For simplicity, we pick $\phi\in M_A$ to be the evaluation at $w=0$. Then
\[\omega_A(z)^*(\phi)=\phi(\omega_A(z))=\frac{dz_0}{z_0}.\]

On the other hand, the topology of $P^c(A)$ is not hard to determine directly. First of all, it is easy to see that 
$P^c(A)\subset \{z_0\neq 0\}$. Hence $P^c(A)={\cc}^{\times}\times p^c(A)$. Under the affine coordinate 
$\xi_j=z_j/z_0,\ 1\leq j\leq n$,
\[p^c(A)=\{\xi=(\xi_1,\ ...,\ \xi_n)\in {\cc}^n:\  1+\sum_{j=1}^n \xi_jw^j\neq 0\ \text{on}\ {\overline{\bf D}}\}.\]
Now consider the maps $H_t$ on $p^c(A)$ defined by
\[H_t(\xi)=(t\xi_1,\ t^2\xi_2,\ ...,\ t^n\xi_n),\]
where $t\in [0,\ 1]$. Since $1+\sum_{j=1}^n \xi_jw^j\neq 0$ for every $w\in {\overline{\bf D}}$,
$1+\sum_{j=1}^n \xi_j(tw)^j$ doesn't vanish on ${\overline{\bf D}}$ as well. This shows that $H_t$ maps $p^c(A)$ into $p^c(A)$,
and it furnishes a retraction of $p^c(A)$ to the origin $0$.

Since $P^c(A)={\cc}^{\times}\times p^c(A)$, $H^*_d(P^c(A),\ \cc)$ is isomorphic to $H^*_d({\cc}^{\times},\ \cc)$.
In particular, $H^1_d(P^c(A),\ \cc)=\cc\frac{dz_0}{z_0}$, and $H^q_d(P^c(A),\ \cc)=0$ for $q\geq 2$. So Corollary 4.2
holds for the disk algebra $A({\bf D})$ and the tuple $(1,\ w,\ w^2,\ ...,\ w^n)$. It is not clear whether Corollary 4.2
holds for other tuples in $A({\bf D})$. \\

{\bx}

It is worth mentioning that for every $\phi$ in the dual of $A({\bf D})$ such that
$\phi(1)=0$, $\phi(\omega_A(z))$ is a closed form on $p^c(A)$ (by the remarks at the end of section 3). 
Using the affine coordinate $\xi$,
$\phi(\omega_A(z))$ is equal to 
\[\phi^*(\xi):=\phi\left(\frac{\sum_{j=1}^n w^jd\xi_j}{1+\sum_{j=1}^n \xi_jw^j}\right)=d\phi\left(log(1+\sum_{j=1}^n \xi_jw^j)\right).\]
Now since in this case $p^c(A)$ is contractible to a point, $\phi^*(\xi)$ is exact, i.e. $\phi\left(log(1+\sum_{j=1}^n \xi_jw^j)\right)$ is 
holomorphic on $p^c(A)$. 

We end this paper with the following\\
 
{\bf Question B}. For the disk algebra $A({\bf D})$, is the map 
\[\omega_A(z)^*:\ \wedge (A^*({\bf D})) \longrightarrow H^*_d(U,\ \cc)\] surjective 
for every tuple $A$ of functions in $A({\bf D})$ and any connected component $U$ of $P^c(A)$?


\vspace{1cm}


\begin{thebibliography}{Abcd}

\large
\bibitem[Ar]{} V. I. Arnold, {\em The cohomology ring of the colored braid group}, 
Mat. Zametki {\bf 5} (1969) 227-231: Math. Notes {\bf 5} (1969) 138-140.\\

\bibitem[Br]{} E. Brieskorn, {\em Sur les groupes de tresses. In: Seminaire Bourbaki 1971/72.}
Lecture Notes in Math. {\bf 317}, Springer Verlag, 1973, 21-44.\\

\bibitem[Da]{} K. R. Davidson, {\em $C^*$-algebras by example}, Fields Institute Monographs,
A.M.S, 1996.\\

\bibitem[Di]{} A. Dimca, {\em Singularities and topology of hypersurfaces},
Universitext, Springer-Verlag, New York, 1992.\\

\bibitem[Do]{} R. G. Douglas, {\em Banach algebra techniques in operator theory},
Pure and Applied Mathematics, Vol. 49. Academic Press, New York-London, 1972.\\

\bibitem[Gr]{} P. Griffiths, {\em On the periods of certain rational integrals I, II},
Ann. Math. {\bf 90}(1969), 460-541.\\

\bibitem[H\"{o}]{} L. H\"{o}rmander, {\em An introduction to complex analysis in several variables}, 3rd ed.,
North Holland, Amsterdam, 1990.\\

\bibitem[OT]{} P. Orlik and H. Terao, {\em Arrangements of hyperplanes}, Grundlehren  der mathematischen Wissenschaften 300, Springer-Verlag Berlin 
Heidelberg 1992.\\

\bibitem[Ra]{} R. M. Range, {\em Holomorphic functions and integral representations in
several complex variables}, Graduate Texts in Mathematics,
Springer-Verlag, New York, 1986.\\

\bibitem[Ta]{} J. L. Taylor, {\em A joint spectrum for several commutative operators}, J. Functional Analysis 6 1970 172--191.\\

\bibitem[Ya]{} R. Yang, {\em Functional spectrum of contractions}, J. Funct. Anal. {\bf 250} (2007), 68-85.\\
\end{thebibliography}
\end{document}